\documentclass[10 pt, conference]{ieeeconf}
\IEEEoverridecommandlockouts                              

\usepackage{cite}
\expandafter\let\csname proof\endcsname\relax
\expandafter\let\csname endproof\endcsname\relax
\usepackage{amsmath,amssymb,amsfonts}
\usepackage{graphicx}
\usepackage{textcomp}
\usepackage{xcolor}
\usepackage{amsthm}
\usepackage{subfigure}
\usepackage{hyperref}
\usepackage[ruled]{algorithm2e} 
\usepackage{algorithmic}
\usepackage{booktabs}
\usepackage{balance}
\def\BibTeX{{\rm B\kern-.05em{\sc i\kern-.025em b}\kern-.08em
    T\kern-.1667em\lower.7ex\hbox{E}\kern-.125emX}}
    \newtheorem{theorem}{Theorem}

\newtheorem{lemma}{Lemma}

\newtheorem*{remark}{Remark}
\newtheorem{definition}{Definition}
\newtheorem{proposition}{Proposition}

\usepackage{threeparttable}

\title{\LARGE \bf
Safety Verification and Controller Synthesis for Systems with Input Constraints
}
\author{Han Wang, Kostas Margellos, Antonis Papachristodoulou
\thanks{The authors are with the Department of Engineering Science, University of Oxford, Oxford, United Kingdom. E-mails: {\tt\small \{han.wang, kostas.margellos, antonis\}@eng.ox.ac.uk}}
}

\begin{document}

\maketitle
\thispagestyle{empty}
\pagestyle{empty}
\begin{abstract}
In this paper we consider the safety verification and safe controller synthesis problems for nonlinear control systems. The Control Barrier Certificates (CBC) approach is proposed as an extension to the Barrier certificates approach. Our approach can be used to characterize the control invariance of a given set in terms of safety of a general nonlinear control system subject to input constraints. From the point of view of controller design, the proposed method provides an approach to synthesize a safe control law that guarantees that the trajectories of the system starting from a given initial set do not enter an unsafe set. Unlike the related control Barrier functions approach, our formulation only considers the vector field within the tangent cone of the zero level set defined by the certificates, and is shown to be less conservative by means of numerical evidence. For polynomial systems with semi-algebraic initial and safe sets, CBCs and safe control laws can be synthesized using sum-of-squares decomposition and semi-definite programming. Examples demonstrate our method. 
\end{abstract}

\section{Introduction}\label{sec:introduction}

Safety-critical systems are commonly used in modern autonomous applications, such as unmanned aerial vehicles, autonomous driving and surgical robotics \cite{guiochet2017safety}. Their safety-critical nature requires the behaviour of these systems to remain within a given safe set for an infinite time horizon. For a model of these systems, such a property is straightforwardly related to reachability analysis and reach-avoid games \cite{lygeros2004reachability}, \cite{margellos2011hamilton}, i.e. finding an initial set so that trajectories reach a target set without entering an unsafe region. However, verify safety for general nonlinear systems using these methods is hard due to the computational difficulty of solving the underlying Hamilton Jacobi Isaacs PDE, especially when control actuation constraints are considered.

To overcome this issue, the connection between forward invariance and safety was established in \cite{prajna2005necessity}. Forward invariance is a system-set property which guarantees that the trajectories entering a set cannot escape it \cite{blanchini1999set}. By finding an invariant subset of a safe region, the system is ensured to be safe. To identify a candidate invariant set, the Barrier certificates approach which takes the invariant set as the certificate's sub-zero level set, was proposed in \cite{prajna2007framework}, \cite{prajna2004safety}. Although the properties of this framework have been demonstrated for autonomous systems with and without stochasticity, there is no systematic formulation for the case where control inputs are present. To address this issue, the control Barrier functions (CBF) approach was proposed in \cite{ames2016control}.

Control Barrier functions are a class of functions that are negative in the unsafe regions, and can be used to verify the safety property. Unlike Lyapunov-like Barrier certificates, control Barrier functions are less restrictive by introducing an additional relaxation term in the constraint. Forward invariance is proved by satisfying the constraints and utilizing the comparison lemma \cite{vidyasagar2002nonlinear}. The approach can be easily combined with the control Lyapunov functions approach \cite{freeman1996control} under a unified quadratic programming framework that compromises safety and controller performance \cite{ames2016control}. It was also shown to be applicable and promising in many applications such as adaptive cruise control \cite{xu2017correctness}, bipedal robots \cite{hsu2015control}, multi-robot collision avoidance \cite{chen2017obstacle} and others. 

Later on, extensive methods to improve the feasibility when input limits are taken into account were proposed, such as adaptive CBF \cite{xiao2021adaptive}, \cite{zeng2021safety}, higher relative degree CBF \cite{xiao2019control}, backup CBF \cite{chen2021backup}, singular CBF \cite{tan2021high}. These methods aim at addressing the cases where the CBF based QP is infeasible. Many times a CBF is assumed to be constructed directly from a physical property such as kinodynamics of the vehicle. How to synthesize a CBF efficiently is still an open question, and has attracted significant attention in recent years.

Direct numerical synthesis approaches by sum-of-squares programming \cite{wang2018permissive}, \cite{xu2017correctness}, machine learning \cite{srinivasan2020synthesis}, and deep learning \cite{robey2020learning} have been proposed. All these methods, either via convex optimisation, or learning techniques, consider the synthesis of a CBF with a relaxation term included in the synthesis procedure. From the standpoint of control invariant sets, it is guaranteed that there exists a class-$\mathcal{K}$ relaxation term to bound the safety variation, but imposing such a term at every point inside the set during the control synthesis introduces conservativeness. Abandoned this term during the synthesis process has been considered in \cite{clark2021verification}, and using Positivstellensatz, a weaker condition on invariance is imposed for systems without input limits. 

In this work we revisit the Barrier certificates approach, and extend it for nonlinear control systems with actuation constraints. Our formulation is a direct interpretation of control invariance and safety guarantee, thus alleviating conservativeness. The existence of a CBC is proved to be sufficient to guarantee safety, hence the approach can be used for safety verification. For systems with polynomial dynamics and semi-algebraic safe and  initial sets, we use sum-of-squares programming and the generalised S-procedure to synthesize a CBC, as well as a Lipschitz continuous safe control law which fulfills the actuation constraints.

The remainder of this paper is organized as follows. The notion of control Barrier certificates is introduced in Section \ref{sec:CBC}. The computation methods with sum-of-squares programming and the S-procedure is presented in Section \ref{sec:computation}. Several simulation results on synthesizing CBCs and safe controllers are shown in Section \ref{sec:simulation}. Section \ref{sec:conclusion} concludes the paper.

\section{Control Barrier Certificates}
\label{sec:CBC}
In this section, we consider the controller synthesis problem under the promise of safety for nonlinear systems. Existing work on Barrier certificates synthesis either limits the analysis to noisy autonomous systems, or tries to design a control law in an online quadratic programming framework. There is no work focusing on combining Barrier certificates construction with controller design, which only requires safety on the boundary of the invariant set. Here, we extend the results of Barrier certificates to control Barrier certificates for safety verification and safe controller design. We also compare our results with the CBF approach. 

\textit{Notation}: $\mathbb{R}$ represents the space of real numbers, and $\mathbb{R}^n$ denotes the $n-$dimensional real space. For a set $S$, $\mathrm{Int}S$, $\partial S$, $\bar S$ are the interior, boundary and complementary set, respectively. $A\succeq 0$ means matrix $A$ is positive semi-definite. $\Sigma[x]$ and $\mathcal{R}[x]$ denote the set of sum-of-squares polynomials and polynomials in $x$ with real coefficients.

\subsection{Control Barrier Certificates Formulation}
We start the formulation for a continuous-time nonlinear system for generality. The system is described by an ordinary differential equation:
\begin{equation}\label{eq:consys}
    \dot x=f(x,u),
\end{equation}
where $x(t) \in \mathbb{R}^n$ denotes the state vector, and $u(t) \in \mathcal{U} \subseteq \mathbb{R}^m$ is the control input, where $\mathcal{U}$ is a bounded set denoting actuation limits and $f(\cdot,\cdot)$ is a locally Lipschitz continuous vector field. We assume that the solution to \eqref{eq:consys} is unique. The flow $\psi(x,t,u)$ denotes the solution of (\ref{eq:consys}) at time $t$ from initial condition $x$ under control $u$. The definitions of a reachable set, forward invariance and safety can be extended to the control system setting.
\begin{definition}[Control Reachable Set]\label{def:conreach}
Consider a vector field $f(\cdot,\cdot)$, a set $X\subseteq \mathbb{R}^n$ and time horizon $T\in\mathbb{R}$. Then the control reachable set of $X$ with respect to vector field $f(\cdot,\cdot)$, control law $u$ and time horizon $T$ is $R_{f,u}^T(X):=\{\psi(x,t,u)|x\in X, t\in T, u\in \mathcal{U}\}.$
\end{definition}
We note here that although the vector field $f(\cdot,\cdot)$ in \eqref{eq:consys} already includes the control input $u$, we still denote the input explicitly in the subscript to distinguish this from the reachable set $R_{f}^T(X)$ for the autonomous system $\dot x=f(x)$.

\begin{definition}[Control Invariant Set]\label{def:coninvariant} A set $X$ is said to be control invariant with respect to vector field $f(\cdot,\cdot)$ if there exits $u$, such that $R_{f,u}^{\infty}(X)\subseteq X$.
\end{definition}
If $u=0$, we call the set $X$ \emph{positive invariant}. The difference between positive invariance and control invariance is obvious: the control effort allows guaranteeing that the flow stays in the set. Hence, the safety of the control system not only depends on the vector field and the predefined safe set $S$, but also on the control admissible set $\mathcal{U}$. 

\begin{definition}[Safety]\label{def:consafety}
Given system \eqref{eq:consys}, an initial set $I$ and a safe set $S$, we say that the system is \emph{safe} if there exits $u\in \mathcal{U}$ such that $R_{f,u}^{\infty}(I)\cap \bar S=\emptyset$.
\end{definition}

The definition of safety of a controlled system is similar to that of an autonomous system. To incorporate safety for the nonlinear control system \eqref{eq:consys}, we aim at finding controller $u$ and a control invariant set $W$, which includes $I$ in its interior and is a subset of the safe set $S$. In particular, $W$ and $u$ fulfill:
\begin{subequations}\label{eq:conalterinvariant}
\begin{align}
    &R_{f,u}^{\infty}(W)\subseteq W,\label{eq:conalterinvariant:invariant}\\
    &R_{f,u}^{\infty}(I)\subseteq W,\label{eq:conalterinvariant:subset}\\
   &W \subseteq S.\label{eq:conalterinvariant:disjoint}
\end{align}
\end{subequations}

\begin{lemma}\label{lem:consubsetinvariant}
If there exits a set $W$ and control input $u\in\mathcal{U}$, such that conditions \eqref{eq:conalterinvariant:invariant}-\eqref{eq:conalterinvariant:disjoint} hold for a vector field $f(x,u)$, then \eqref{eq:consys} is safe according to Definition \ref{def:consafety}.
\end{lemma}

We now define the notion of \emph{Control Barrier Certificates} (CBC) for finding a feasible candidate control invariant set, and a controller according to condition \eqref{eq:conalterinvariant}.

\begin{definition}[Control Barrier Certificates]\label{Def: CBC}
Let a continuous time control system denoted by $\dot x=f(x,u)$, with initial set $I\subseteq \mathbb{R}^n$, safe set $S\subseteq \mathbb{R}^n$, and input constraints $\mathcal{U}\subseteq \mathbb{R}^m$. A $C^1$ function $B:\mathbb{R}^n\to \mathbb{R}$ is called a \emph{Control Barrier Certificate (CBC)} if 
\begin{subequations}\label{eq:CBC}
\begin{equation}\label{eq:CBC1}
    B(x)<0,~~~~~\forall x\in \bar S,
\end{equation}
\begin{equation}\label{eq:CBC2}
    B(x)\ge0,~~~~~\forall x\in I,
\end{equation}
\begin{equation}\label{eq:CBC3}
    \sup_{u\in\mathcal{U}}\frac{\partial B(x)}{\partial x}f(x,u)>0,~~~~\forall x\in \partial{B}.
\end{equation}
\end{subequations}
\end{definition}
\vspace{1ex}

Let 
\begin{equation}\label{eq:KCBC}
    {K_{CBC}}(x): = \left\{ \begin{array}{l}
\left \{ u|\frac{{\partial B(x)}}{{\partial x}}f(x,u) > 0\ \right \}  \cap {\cal U},~\mathrm{if}~B(x)=0\\
{\cal U,~\mathrm{otherwise}}
\end{array} \right.
\end{equation}
denote the admissible set of control inputs for a CBC $B(x)$. Let $\mathcal{B}:=\{x|B(x)\ge0\}$ denote the  zero-super level set of $B(x)$. We then have the following result on safety.

\begin{theorem}\label{th:CBC}
Consider \eqref{eq:consys}, a safe set $S$ and an initial set $I$. If there exists a CBC $B(x)$ that satisfies conditions \eqref{eq:CBC}, then for any state $x$ and any $u\in K_{CBC}(x)$, the safety of system \eqref{eq:consys} is guaranteed.
\end{theorem}

\begin{proof}
Equation \eqref{eq:CBC1} indicates that for any $x\in\bar{S}$, we have $x\in\bar{\mathcal{B}}$, thus $\mathcal{B}\subseteq S$, which shows that condition \eqref{eq:conalterinvariant:disjoint} holds. Similarly, Equation \eqref{eq:CBC2} demonstrates that $I\in\mathcal{B}$. Therefore, we only need to prove that $R_{f,u}^{\infty}(\mathcal{B})\in\mathcal{B}$ to show conditions  \eqref{eq:conalterinvariant:invariant} -- \eqref{eq:conalterinvariant:subset}. We recall that under control input $u$, the vector field $f(x,u)$ is locally Lipchitz continuous, and the solution is unique. This indicates that the flow $\psi(x,t,u)$ is continuous over $t$. Besides, the fact that $B(x)$ is a $C^1$ function guarantees that trajectories starting from $\mathrm{Int}(\mathcal{B})$ to $\bar {\cal B}$ will cross $\partial \mathcal{B}$. Thus, any bounded input $u\in\mathcal{U}$ at $x\in \mathrm{Int}(\mathcal{B})=\{x|B(x)>0\}$ shows the positivity of $B(\psi(x,t,u))$ when $t\to 0$. Regarding the boundary, for any $x\in \partial \mathcal{B}$, $\dot B(x)=\frac{\partial B(x)}{\partial x}\frac{dx}{dt}=\frac{\partial B(x)}{\partial x}f(x,u)>0$ from the definition of CBC. Thus, the vector field $f(x,u)\in \mathrm{Tang}_\mathcal{B}(x)$ for $x\in\partial \mathcal{B}$ and $u\in\mathcal{U}$. According to the subtangenality Theorem \cite{nagumo1942uber}, the set $\mathcal{B}$ is control invariant with vector field $f(x,u)$, which directly indicates that $\mathcal{B}$ is control invariant. According to Lemma \ref{lem:consubsetinvariant}, $\mathcal{B}$ is a feasible candidate control invariant set verifying the safety of the control system \eqref{eq:consys}.
\end{proof}

Theorem \ref{th:CBC} shows that the existence of a control Barrier certificate $B(x)$ ensures safety for safety-critical systems. Meanwhile, the control admissible set \eqref{eq:KCBC} certifies the selection of control effort. For problems that the control Barrier certificate $B(x)$ can be easily constructed and verified through physical properties, one can formulate a quadratic program to synthesize the safe controller $u$ at $x$.
\begin{equation}\label{eq:CBCQP1}
    \begin{split}
        &\mathop{\min }\limits_u ||u - {u^*(x)}||\\
    &\mathrm{s.t.}~u\in K_{CBC}(x),      
    \end{split}
\end{equation}
where $u^*(x)$ is a nominal control input designed from other tools, e.g. PID, MPC. We note here the formulation \eqref{eq:CBCQP1} is different from that of CBF based QP. Here, in the interior of the control invariant set $\mathcal{B}$, the solution of \eqref{eq:CBCQP1} is the direct projection from $u^*(x)$ on the control admissible set $\mathcal{U}$. 

For the scenario where the control Barrier certificate is unknown, the problem is to synthesise the control Barrier certificates together with the safe controller design. To begin with, $B(x)$ is parameterized by 
\begin{equation}\label{eq:parameterize}
    B(x)=p_1\Lambda_1(x)+\ldots+p_k\Lambda_k(x),
\end{equation}
where $p:=\{p_1,\ldots,p_k\}$ denotes a series of parameters which will be decision variables in an optimisation problem, and $\Lambda_1(x),\ldots,\Lambda_k(x)$ are a class of function basis. The new optimisation problem for constructing the CBC and controller is
\begin{equation}\label{eq:CBCQP2}
\begin{split}
    \mathrm{find}~&u(x),p\\
    \mathrm{s.t.}~&\eqref{eq:parameterize},\eqref{eq:CBC};\\
    &u(x)\in K_{CBC}(x).
\end{split}
\end{equation}
Compared to quadratic programming \eqref{eq:CBCQP1} with known control Barrier certificates, \eqref{eq:CBCQP2} is computationally intractable since it involves solving an infinitely constrained optimisation problem. We will show how to address this difficulty by sum-of-squares programming in Section \ref{sec:computation}.

\section{Computation Method}
\label{sec:computation}
In this section we show how to construct the control Barrier certificates and the safe control law for polynomial systems with semi-algebraic safe and initial sets. The nonlinear control affine system is represented by
\begin{equation}\label{eq:nonlinaff}
   \dot x=f(x)+g(x)u, 
\end{equation}
where $f(x)$ and $g(x)$ are locally smooth polynomial functions and $u\in \mathcal{U}:=\{u|Au+b\ge 0\}$.

Even for such a simplified system model, solving the parametric optimisation problem \eqref{eq:parameterize} -- \eqref{eq:CBCQP2} involves solving an infinite set of non-negative inequalities and hence is computationally intractable. However, for systems with polynomial functions $f(x)$, $g(x)$ and semi-algebraic sets $I$, $S$, a tractable method for tackling the infinite inequalities is \textit{sum-of-squares} (SOS) programming, which is a convex relaxation method based on the sum-of-squares decomposition of multivariate polynomials and semidefinite programming. 

A SOS program is a convex optimisation problem of the following form:
\begin{equation}\label{eq:SOS}
 \begin{split}
     &\mathop{\min} \limits_p \sum^k_{j=1}w_jp_j\\
     &\mathrm{s.t.}~h_0(x)+\sum^k_{j=1}p_jh_j(x)\in\Sigma[x],
     \end{split}   
\end{equation}
where the decision variables $p_1,\ldots,p_k$ are real parameters, and $w_1,\ldots,w_k$ are predefined weight constants. Also, [$h_0(x),\ldots,h_k(x)$] is a polynomial basis in $x$. A multivariate polynomial $\mathrm{s.t.}~h_0(x)+\sum^k_{j=1}p_jh_j(x)$ with $x\in\mathbb{R}^n$ is a SOS polynomial if there exists $k$ polynomials $f_1(x)\ldots f_k(x)$ such that $f(x)=\sum_{i=1}^kf_i^2(x)$. Then it directly follows that a SOS $f(x)$ is non-negative for any $x\in\mathbb{R}^n$. A SOS program can be transformed into a semi-definite program with $f(x)=Z^\top QZ(x)$, where $Q\succeq 0$ and $Z(x)$ is a monomial vector. 

\subsection{SOS for CBC Synthesis}

To interpret the constraints \eqref{eq:CBC} into SOS constraints, we assume that the resulting control Barrier certificate is a polynomial function parameterized by real coefficients $p_1,\ldots,p_m$ in the following way
\begin{equation}\label{eq:paraCBC}
    B(x)=p_0+\sum_{j=1}^mp_jb_j(x),
\end{equation}
where $b_j(x)$s are polynomial or monomial function bases, and $p_0$ is a positive real scalar. Similarly, the control input is parameterized by real scalar coefficients $k_1,\ldots, k_l$, and a real vector coefficient $k_0 \in \mathbb{R}^m$ with
\begin{equation}\label{eq:paracontrol}
    u(x)=k_0+\sum_{j=1}^lk_jv_j(x),
\end{equation}
where $v_j(x)$s are polynomial or monomial vector bases. We note here there the reason why we use the constant term $k_0$ is different from that of $p_0$. From the view of control, $k_0$ introduces a feedforward term, which in some cases is important for safety, for example at some singular points where $\mathop{\sum \limits_{j=1}^lk_jv_j(x)}=0$.

\begin{theorem}\label{th:SOSCBC}
Consider a polynomial nonlinear system \eqref{eq:nonlinaff}, semi-algebraic safe set $S=\{x|s(x)\ge 0\}$, initial set $I=\{x|w(x)\ge 0\}$, and control admissible set $\mathcal{U}:=\{u|Au+b\ge0\}$, where $A\in \mathbb{R}^{h\times m}$, and $b\in \mathbb{R}^h$. If there exit multipliers $\sigma_{\mathrm{safe}} \in \Sigma [x]$, $\sigma_{\mathrm{init}}\in \Sigma[x]$, $\lambda_1 \in \mathcal{R}[x]$, $\lambda_2\in \mathcal{R}[x]^h$, polynomials $B(x)\in\mathcal{R}[x]$, $u(x)\in\mathcal{R}[x]$, predefined small positive real scalars $\epsilon_1>0$, $\epsilon_2>0$, such that 
\begin{subequations}\label{eq:SOSCBC}
\begin{align}
    &-B(x)+\sigma_{\mathrm{safe}}s(x)-\epsilon_1 \in \Sigma[x],\label{eq:SOSCBC1}\\
    &B(x)-\sigma_{\mathrm{init}}w(x)\in \Sigma[x],\label{eq:SOSCBC2}\\
    &\frac{\partial B(x)}{\partial x}(f(x)+g(x)u(x))+\lambda_1 B(x)-\epsilon_2 \in \Sigma[x],\label{eq:SOSCBC3}\\
    &-\lambda_2 B(x)+Au(x)+b \in \Sigma[x]^h,\label{eq:SOSCBC4}
\end{align}
\end{subequations}
then $B(x)$ fulfills the conditions $\eqref{eq:CBC}$ and $\mathcal{B}=\{x|B(x)\ge 0\}$ is a control invariant set with respect to vector field $f(x)+g(x)u(x)$. 
\end{theorem}

\begin{proof}
Condition \eqref{eq:SOSCBC1} indicates that for any $x$, $-B(x)+\sigma_{\mathrm{safe}}s(x)-\epsilon_1\ge0$, thus for any $x$, $-B(x)+\sigma_{\mathrm{safe}}s(x)>0$. Therefore, for any $x\in \bar{S}$, we directly have that $\sigma_{\mathrm{safe}}s(x)\le0$, and further $B(x)< 0$, i.e., \eqref{eq:CBC1} holds. Similarly \eqref{eq:SOSCBC2} can be shown to satisfy \eqref{eq:CBC2} following the same arguments. Based on the S-procedure, condition \eqref{eq:SOSCBC3} implies condition \eqref{eq:CBC3}, because when $B(x)=0$, $\frac{\partial B(x)}{\partial x}(f(x)+g(x)u(x)) -\epsilon_2 \ge 0$, and thus $\frac{B(x)}{\partial x}(f(x)+g(x)u(x))>0$. Condition \eqref{eq:SOSCBC4} implies that $Au(x)+b$ is elementary-wise nonnegative for $x\in \partial \mathcal{B}$. The small positive real scalars $\epsilon_1$, $\epsilon_2$ ensure strict inequality for \eqref{eq:CBC1} and \eqref{eq:CBC3}.
\end{proof}

We note that in Theorem \ref{th:SOSCBC} we only require a polynomial multiplier $\lambda$, but not a SOS one since the condition $\frac{\partial B(x)}{\partial x}(f(x)+g(x)u(x))\ge 0$ is only imposed on the boundary $B(x)=0$. Condition \eqref{eq:SOSCBC3} introduces products of decision variables, i.e. $\lambda B(x)$, which results in bilinearity. However, there is no guaranteed solver for nonconvex, or specifically bilinear constrained SOS programs. Here, like existing work of using SOS to synthesize Barrier certificates, we use an iterative procedure for control Barrier certificate synthesis and safe control law design. Different from the iterative algorithm for Barrier certificate synthesis, our problem involves an additional polynomial variable $u$ in the SOS program. Thus, an additional round for controller synthesis is required in our algorithm.

1) \textbf{Initialization}: We first fix the degree of polynomials $B(x)$, $\sigma_{\mathrm{safe}}$, $\sigma_{\mathrm{init}}$, $\lambda_1$, $\lambda_2$ and $u(x)$. The polynomial/monomial scalar/vector bases $b_j(x)$s and $v_j(x)$s in \eqref{eq:paraCBC} and \eqref{eq:paracontrol} have degree upper bounded by the aforementioned degrees of $B(x)$. $\epsilon_1$ and $\epsilon_2$ are chosen to be small real numbers. Unlike the iterative procedure proposed in \cite{xu2017correctness} which initializes the control law by a scaled LQR controller, we find the initialized feasible control input $u^0(x)$ by solving a feasibility SOS program.
\begin{equation}\label{eq:InitSOSCON}
\begin{split}
     &\mathrm{find~}k_1,\ldots, k_l,\sigma_{\mathrm{cont}}\\
      &\mathrm{s.t.}~A(k_0+\sum_{j=1}^lk_jv_j(x))+b\cdot \sigma_{\mathrm{cont}}\in\Sigma[x]^h.
\end{split}
\end{equation}
We note here that there is no assumed control Barrier certificate $B(x)$ at this stage of finding the initial feasible control input $u^0(x)$. Therefore, $u^0(x)$ can not be restricted to the domain of $\partial \mathcal{B}$ as that in \eqref{eq:SOSCBC4}. Other than directly interpreting $A(k_0+\sum_{j=1}^lk_jv_j)+b \in \Sigma[x]^h$, we add an additional positive multiplier $\sigma_{\mathrm{cont}}$ which satisfies $\sigma_{\mathrm{cont}}-\epsilon_3\in \Sigma[x]$, $\epsilon_3>0$ to avoid introducing constant terms in the SOS constraints, as well as improving feasibility. The resulting initial controller $u^0(x)$ is derived by the parameters $k_1,\ldots,k_l$ and the scaled term $\sigma_{\mathrm{cont}}$ from the solution of \eqref{eq:InitSOSCON}
\begin{equation}\label{eq:initcon}
    u^0(x)=\frac{1}{\sigma_{\mathrm{cont}}}\cdot(k_0+\sum_{j=1}^lk_jv_j(x)).
\end{equation}
The feasibility of such an initialized controller is guaranteed by the following proposition.
\begin{proposition}\label {pro:initcon}
The initialized control input $u^0(x)$ in \eqref{eq:initcon} satisfies $Au^0(x)+b\ge0$.
\end{proposition}
\begin{proof}
We have $A(k_0+\sum_{j=1}^lk_jv_j)+b\cdot \sigma_{\mathrm{cont}}\ge 0$ from the SOS constraints $A(k_0+\sum_{j=1}^lk_jv_j)+b\cdot \sigma_{\mathrm{cont}} \in \Sigma[x]^h$ in \eqref{eq:InitSOSCON}. Because of the positivity of the multiplier $\sigma_{\mathrm{cont}}$, we directly have $A(\frac{1}{\sigma_{\mathrm{cont}}}\cdot(k_0+\sum_{j=1}^lk_jv_j))+b\ge 0$, which indicates $u^0(x)\in\mathcal{U}$.
\end{proof}

Given initial input $u^0(x)$, the corresponding scaled multiplier $\sigma_{\mathrm{cont}}$, the initial control Barrier certificate $B^0(x)$ can be found by solving an initial feasibility SOS program as
\begin{equation}\label{eq:InitSOSCBC}
    \begin{split}
     \mathrm{find}~&p_0,\ldots ,p_m,\sigma_{\mathrm{safe}},\sigma_{\mathrm{init}}\\
        \mathrm{s.t.}~ &-B(x)+\sigma_{\mathrm{safe}}s(x)-\epsilon_1\in \Sigma[x],\\
        &B(x)-\sigma_{\mathrm{init}}w(x)\in \Sigma[x],\\
        &\sigma_{\mathrm{cont}}\cdot\frac{\partial B(x)}{\partial x}(f(x)+g(x)u^{0}(x))-\epsilon_2\in \Sigma[x],\\
      &~~~~~~~~~~B(x)~\mathrm{from}~\eqref{eq:paraCBC}.\\
    \end{split}
\end{equation}
The boundary condition \eqref{eq:SOSCBC3} is strengthened to be $\frac{\partial B(x)}{\partial x}(f(x)+g(x)u(x))-\epsilon_2\in\Sigma[x]$ for convexity and simplicity of computing. This condition is also referred to be the weak Barrier certificate in \cite{prajna2004safety}. $\sigma_{\mathrm{cont}}\cdot \frac{\partial B(x)}{\partial x}(f(x)+g(x)u^{0}(x))-\epsilon_2$ is guaranteed to be a polynomial, since $\sigma_{\mathrm{cont}}\cdot u^0(x)$ is a polynomial.

After obtaining a feasible initial control input $u^0(x)$ and control Barrier certificate $B^0(x)$, the problem of control Barrier certificates synthesis can be regarded as a Barrier certificates synthesis problem with vector field $f(x)+g(x)u^0(x)$. The multipliers $\lambda_1^0$, $\lambda_2^0$ are fixed to be 0 or 1 in initialization for simplicity. The initial control Barrier certificate $B^0(x)$ is used to enlarge the size of the control invariant set incrementally. The following steps of the algorithm iteratively solve the SOS program to address the bisecting terms  $\lambda_1 B(x)$ and $\frac{\partial B(x)}{\partial x}(f(x)+g(x)u(x))$ in \eqref{eq:SOSCBC3}.   

2) \textbf{Update the control input} $u^k(x)$:
At iteration $k$, given a control Barrier certificate from \eqref{eq:InitSOSCBC} (when $k=1$) or \eqref{eq:IterSOSCBC} (when $k\ge2$), the controller synthesis is constrained to \eqref{eq:SOSCBC4}. Fixing $B(x)=B^{k-1}(x)$, a convex programming synthesis procedure for $u^k(x)$ is 
\begin{equation}\label{eq:IterSOSCON}
\begin{split}
     &\mathrm{find~} k_0,\ldots, k_l,\lambda_1,\lambda_2\\
      &\mathrm{s.t.}-\lambda_2B^{k-1}(x)+~A(k_0+\sum_{j=1}^lk_jv_j)+b\in\Sigma[x]^h,\\
      &\frac{\partial B^{k-1}(x)}{\partial x}(f(x)+g(x)u(x))+\lambda_1B^{k-1}(x)-\epsilon_2\in\Sigma[x],
\end{split}
\end{equation}
and we have that $u^k(x)=(k_0+\sum_{j=1}^lk_jv_j).$ Here we use $\lambda_1$ other than $\lambda_{1}^{k-1}$ since $B(x)$ has been substituted by $B^k(x)$, thus there is no bilinear term anymore. By limiting the domain of the controller to $\partial \mathcal{B}$, there is no need to have additional multiplier $\sigma_{\mathrm{cont}}$ as that has been used in initial controller design for feasibility.

3) \textbf{Synthesize the control Barrier certificate} $ B^k(x)$:
After obtaining a feasible control input $u^{k-1}(x)$, the synthesis of a control Barrier certificate $B^k(x)$ relies on fixed multipliers $\lambda_1^{k-1}$, $ \lambda_2^{k-1}$ to bypass the bilinear terms. Searching for $B^k(x)$ and the remaining multipliers follows the following SOS program
\begin{equation}\label{eq:IterSOSCBC}
    \begin{split}
        &\mathrm{find}~p_0,\ldots ,p_m,\sigma_{\mathrm{safe}},\sigma_{\mathrm{init}},\sigma_{\mathrm{enl}}\\
        \mathrm{s.t.}~ &-B(x)+\sigma_{\mathrm{safe}}s(x)-\epsilon_1\in \Sigma[x],\\
        &B(x)-\sigma_{\mathrm{init}}w(x)\in \Sigma[x],\\
        & \frac{\partial B(x)}{\partial x}(f(x)+g(x)u^{k}(x))+\lambda_1^{k-1} B(x)-\epsilon_2\in \Sigma[x],\\
        &-\lambda_2^{k-1}B(x)+Au^{k}(x)+b\in\Sigma[x]^h,\\
        &B(x)-\sigma_{\mathrm{enl}}B^{k-1}(x)\in\Sigma[x],\\ 
          &~~~~~~~~~B(x)~\mathrm{in}~\eqref{eq:paraCBC},\\
    \end{split}
\end{equation}
where $\sigma_{\mathrm{enl}}\in\Sigma[x]$.
Here the control law $u^{k-1}(x)$ is substituted for the variable $u$, and the multipliers $\lambda_1$ $\lambda_2$ are substituted by $\lambda_1^{k-1}$ and $\lambda_2^{k-1}$, respectively. We introduce additional constraints $B(x)-\sigma_{\mathrm{enl}}B^{k-1}(x)\in\Sigma[x]$ to enlarge the volume of the control invariant set $\mathcal{B}^{k}$ by enforcing $\mathcal{B}^{k-1}\subseteq \mathcal{B}^k$. A similar technique is also used in \cite{cunis2021viability}.

4) \textbf{Update the multipliers}:
The multiplier $\lambda_1^k$ updates rely on a fixed control Barrier certificate $B^k(x)$ and input $u^k(x)$. Clearly, there is no bilinearity in the control input update procedure \eqref{eq:IterSOSCON}. The multipliers $\lambda_1^k$ and $\lambda_2^k$ are obtained by directly solving it. There is no need to fix $B(x)$ and re-solve the programming problem.
\begin{remark}\label{rem:SOSCBC}
For the case where \eqref{eq:InitSOSCON} or \eqref{eq:InitSOSCBC} is infeasible, there are two options for ensuring feasibility: (i) Increase the degree of the polynomial bases $v_1,\ldots,v_l$, $b_1,\ldots,b_m$; (ii) Re-solve the problem \eqref{eq:InitSOSCON} with an alternative objective function for a different initialization.
\end{remark}

\section{Simulation Results and Discussion}\label{sec:simulation}

In this section we show numerical simulation results on synthesizing control Barrier certificates and safe controllers under different system settings. The SOS toolbox SOSTOOLS \cite{prajna2002introducing}\cite{sostools} is used with version v401 for parsing the SOS programs, while SeDuMi is used for solving the resulting semidefintie program \cite{doi:10.1080/10556789908805766}. We also give a comparison between the CBC proposed in this paper and CBF mainly from the view point of synthesis.

\subsection{Nonlinear Control Affine Systems}
We first consider a general second order polynomial nonlinear control affine system. This system is defined by 
\begin{equation}\label{eq:conaff}
    \begin{bmatrix}
\dot x_1\\
\dot x_2
\end{bmatrix}=
    \begin{bmatrix}
x_2\\
x_1+\frac{1}{3}x_1^3+x_2
\end{bmatrix}+
    \begin{bmatrix}
x_1^2+x_2+1\\
x_2^2+x_1+1
\end{bmatrix}
\begin{bmatrix}
u_1\\
u_2
\end{bmatrix},
\end{equation}
where the control input is box constrained, i.e. $u_1\in[-1.5,1.5]$, $u_2\in[-1.5,1.5]$. The safe set is defined by a disc $S=\{x|x_1^2+x_2^2-3\le 0\}$, and initial set defined by $I=\{x|(x_1-0.4)^2+(x_2-0.4)^2-0.16\le 0\}$. We leverage the control Barrier certificates synthesis procedures \eqref{eq:SOSCBC}
to find a polynomial CBC $B_1(x)$, and compare the results with the CBF synthesis procedure proposed in \cite{xu2017correctness}. To synthesize a candidate CBF $B_2(x)$, an alternative constraint for \eqref{eq:SOSCBC3} is introduced
\begin{equation}\label{eq:cbf}
    \frac{\partial B(x)}{\partial x}(f(x)+g(x)u(x))-\sigma_{\mathrm{cbf}}B(x)+\alpha B(x)-\epsilon_2\in\Sigma[x],
\end{equation}
where the class-$\mathcal{K}$ function is selected to be $\alpha B(x)$ with $\alpha>0$, and $\sigma_{\mathrm{cbf}}\in\Sigma[x]$ is a SOS multiplier. Here we restrict the definition domain for CBF to be $\mathcal{B}_2$. $\epsilon_2$ is set to be the same with that in \eqref{eq:SOSCBC3}. Instead of the feasibility SOS program used for CBC, we set an objective function $\alpha$ which is maximized for CBF as in \cite{xu2017correctness}.

\begin{figure}[h]
    \centering
    \includegraphics[width=0.4\textwidth]{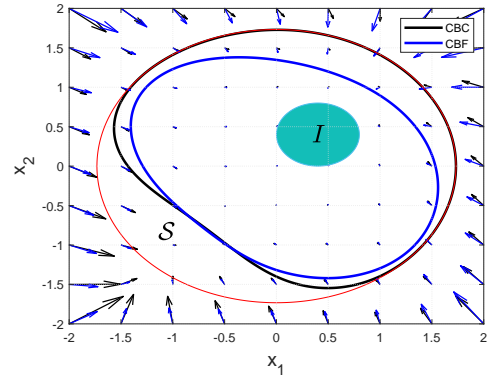}
    \caption{Control invariant sets defined by CBC or CBF}
    \label{fig:comparison}
\end{figure}

Figure \ref{fig:comparison} shows the control invariant sets defined by CBF and CBC. The red and light blue disc represent the safe and initial sets, respectively. The interior of the deep blue curve is the invariant set $\mathcal{B}_2$ defined by CBF, and the interior of the black curve is the invariant set $\mathcal{B}_1$ defined by CBC. It can be seen from the figure that $\mathcal{B}_1$ is ``larger" than $\mathcal{B}_2$. Actually we have $\mathcal{B}_2\subset\mathcal{B}_1$, which is proved by there exists a SOS multiplier $\sigma$, such that $B_1(x)-\sigma B_2(x)\in\Sigma[x]$. The reason is that, we trivially have $\sigma_{\mathrm{cbf}}+\alpha\in\mathcal{R}[x]$. A larger search area enables us to find a larger control invariant set. On the other hand, the additional term $\lambda_1B_1(x)$ can be regarded as an adapted relaxation term compared with a fixed class-$\mathcal{K}$ function used in CBF approach. By using a zeroth order base for the polynomial multiplier $\lambda_1$ and expanding the definition domain of CBF to the whole real space, our formulation is equivalent to CBF. Higher order basis selections hereby reduce conservativeness. 

\begin{figure}[h]
    \centering
    \includegraphics[width=0.4\textwidth]{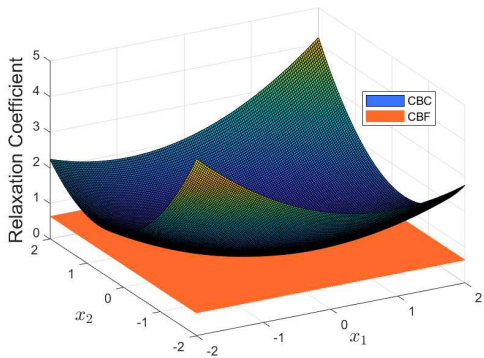}
    \caption{Relaxation coefficients $\lambda(x)$ for CBC, and $\alpha$ for CBF}
    \label{fig:multiplier}
\end{figure}

Figure \ref{fig:multiplier} shows the value of relaxation coefficient $\lambda_1$ and $\alpha$. The multiplier $\lambda_1$ includes the following monomial basis: $[x_1^2,x_1x_2,x_2^2,x_1,x_2,1]$. It can be seen that $\lambda_1$ varies in the control invariant set, which therefore endows the formulation flexibility. An interesting property here is that $\alpha$ cannot be too large, this is because for $x\in\bar{\mathcal{B}_1}$, $\alpha B_2(x)<0$. In addition, with a non-empty safe set $S\subset \mathbb{R}^n$, we directly have $B_1(x)\notin \Sigma[x]$, and $\alpha B_1(x)\notin \Sigma[x]$.

\begin{figure*}[h]
    \centering
        \subfigure[Phase portrait for the system \eqref{eq:conaff} ]{
        \includegraphics[width=0.31\textwidth]{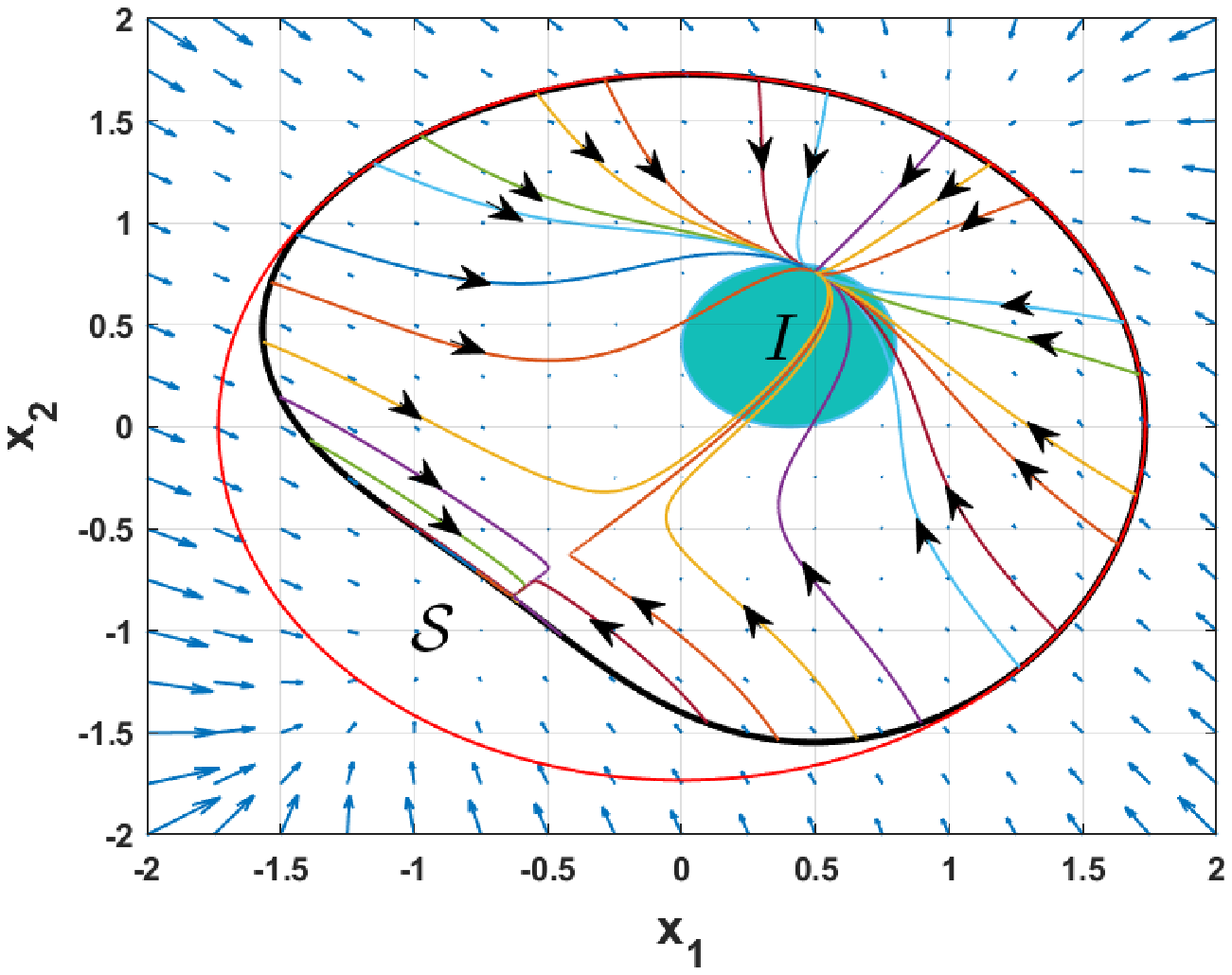}
        \label{fig:Nonlinsys}
    }
    \subfigure[Level set of $u_1$]{
        \includegraphics[width=0.31\textwidth]{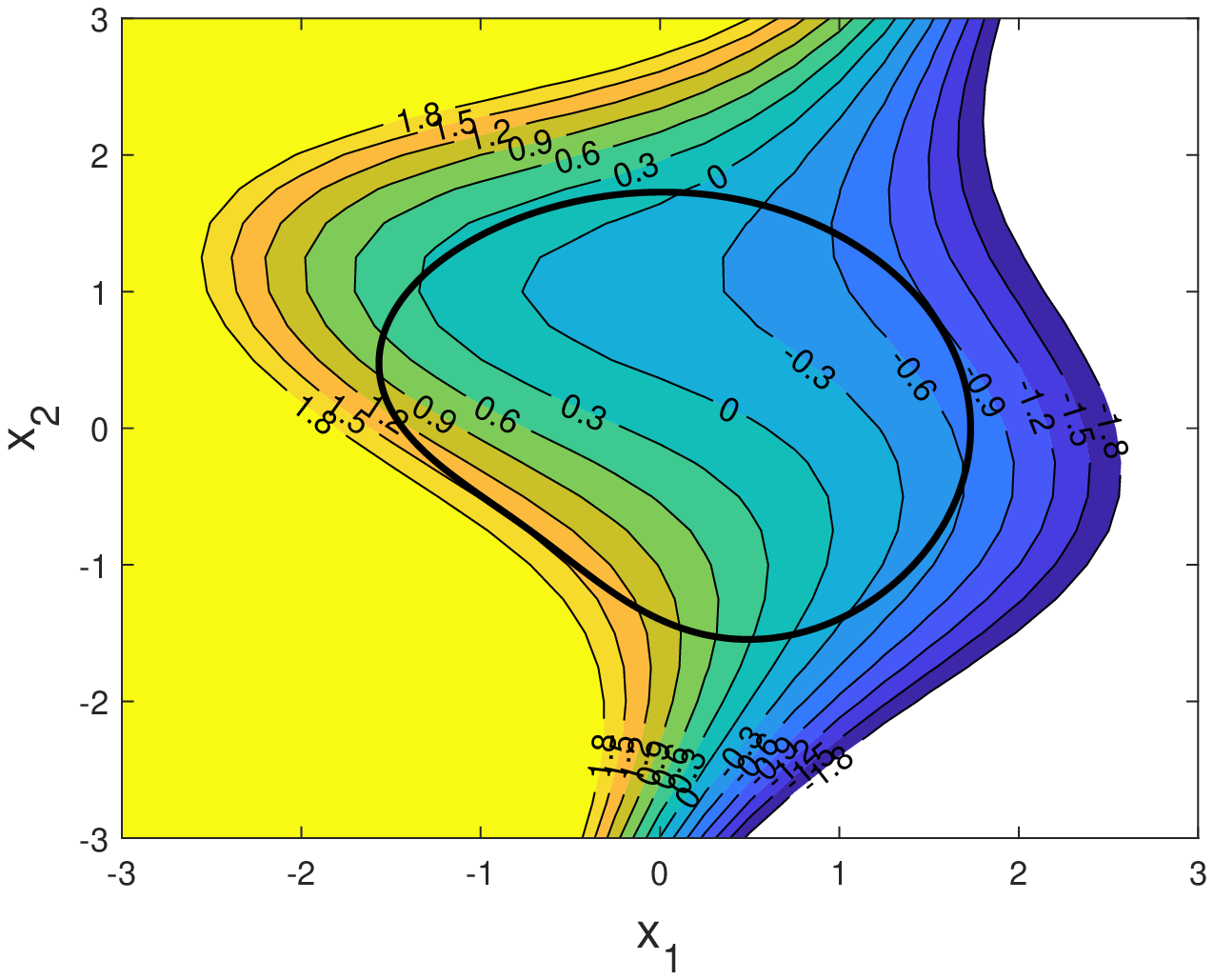}
        \label{fig:NonSolU1}
    }
    \subfigure[Level set of $u_2$]{
        \includegraphics[width=0.31\textwidth]{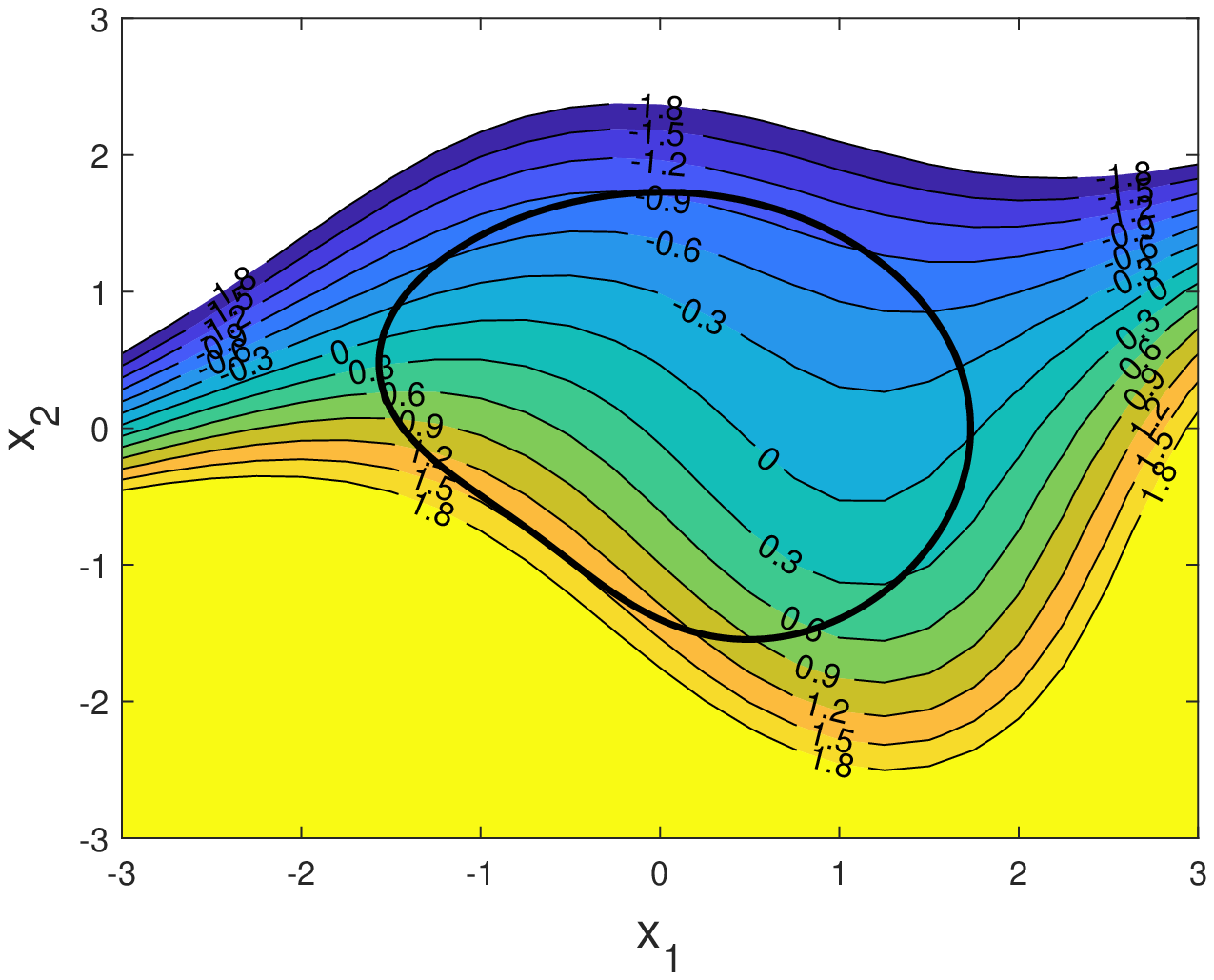}
        \label{fig:NonSolU2}
    }
    \caption{The interior of the red disc represents the safe set, the interior of the blue disc represents the initial set from which the trajectories start. The black closed curve encircling the initial set is the control invariant set, defined by the super-zero level set of $B_1(x)$. The arrows in the figure represent the vector field. The colorful lines are the trajectories starting from $\partial \mathcal{B}_1$.} 
    \label{fig:NonLinSysCon}
\end{figure*}

The control invariant set $\mathcal{B}_1$ obtained by CBC design and values of the safe controllers are shown in Figure \ref{fig:NonLinSysCon}. The vector field, which is represented by the arrows in Figure \ref{fig:Nonlinsys} point inside $\mathcal{B}_1$ on $\partial{\mathcal{B}_1}$. The value of the polynomial control law $u_(x)$ is within $[-1.5,1.5]$ in both coordinates. 

\begin{figure*}[h]
    \centering
        \subfigure[Phase portrait for the system \eqref{eq:linsys} ]{
        \includegraphics[width=0.31\textwidth]{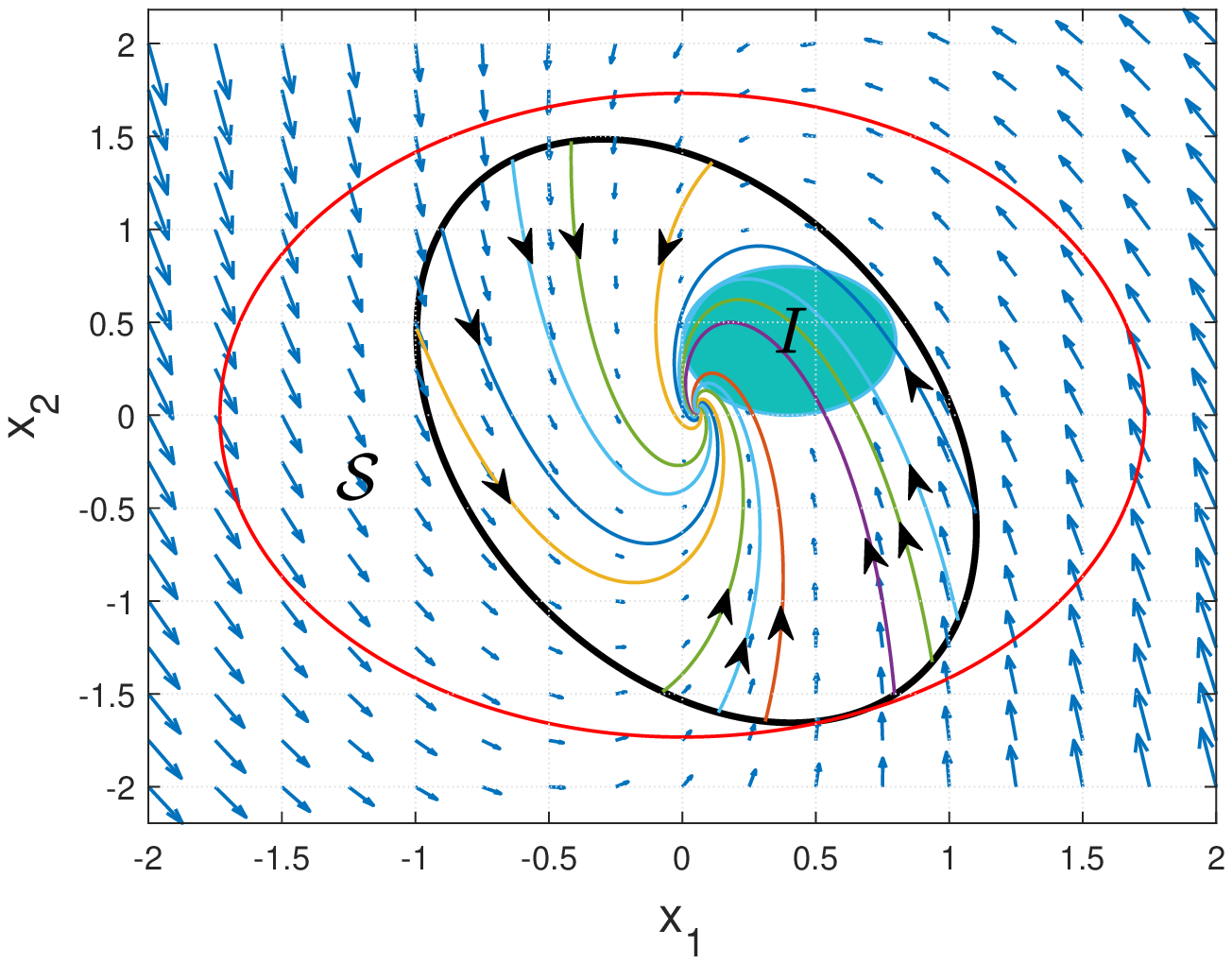}
        \label{fig:linsys}
    }
    \subfigure[Level set of $u_1$]{
        \includegraphics[width=0.31\textwidth]{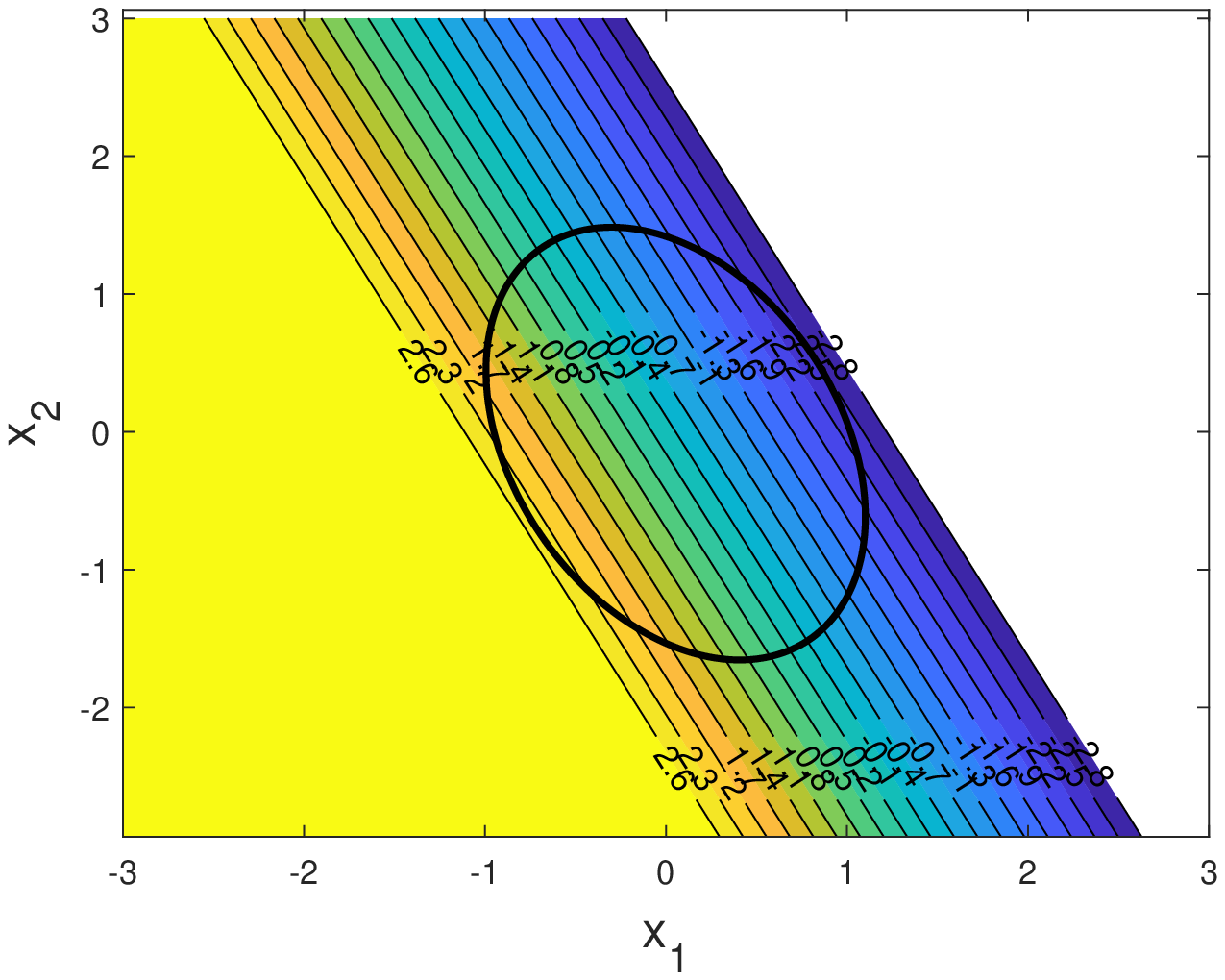}
        \label{fig:SolU1}
    }
    \subfigure[Level set of $u_2$]{
        \includegraphics[width=0.31\textwidth]{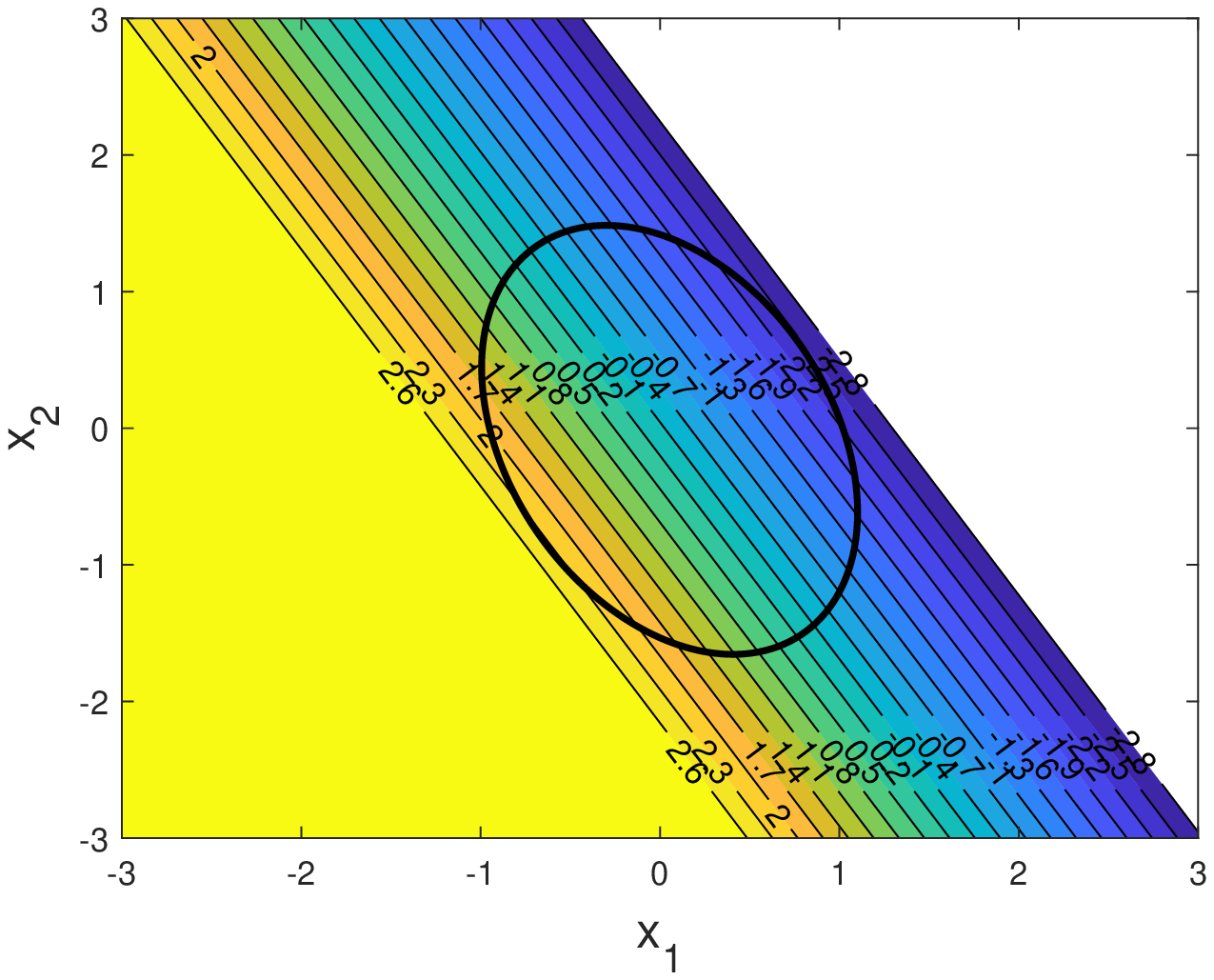}
        \label{fig:SolU2}
    }
    \caption{The safe and initial set are defined to be the same as in Figure \ref{fig:NonLinSysCon}. Safety is ensured with the polynomial control law.}
    \label{fig:LinSysCon}
\end{figure*}

\subsection{LTI Systems}
Consider a second order linear model

\begin{equation}\label{eq:linsys}
      \begin{bmatrix}
\dot x_1\\
\dot x_2
\end{bmatrix}= 
    \begin{bmatrix}
2&1\\
3&1
\end{bmatrix}
    \begin{bmatrix}
x_1\\
x_2
\end{bmatrix}
+
    \begin{bmatrix}
u_1\\
u_2
\end{bmatrix},
\end{equation}
where $u_1\in[-2.5,2.5]$, $u_2\in[-2.5,2.5]$. The system is unstable since the eigenvalues of the state matrix $    \begin{bmatrix}
2&1\\
3&1
\end{bmatrix}$ are $3.3$ and $-0.3$, whereas it is locally stabilizable. The safe set is defined by a disc $S=\{x|x_1^2+x_2^2-3\le 0\}$. The trajectories of the system start from the following initial set $I=\{x| (x_1-0.4)^2+(x_2-0.4)^2-0.16\le 0\}$. Clearly, all trajectories starting from the initial set tend to infinity, since the system is unstable. Safety is therefore violated with a closed safe region set.

Using a second degree basis $[1,x_1,x_2,x_1x_2,x_1^2,x_2^2]$, a feasible candidate CBC is given by $B_1(x)=  -7.635x_1^2 - 3.439x_1x_2 - 3.4024x_2^2 + 0.5x_1 
  - 0.4x_2 +7.402$. The corresponding control inputs lying inside $[-2.5,2.5]$ when $x
\in\partial \mathcal{B}_1$ are $u_1(x)=-2.32x_1-1.11x_2+0.022$, $u_2(x)=-2.12x_1-1.27x_2-0.046$. Obviously $u_1(x)$ and $u_2(x)$ are admissible only within some local regions. More specifically, within $\mathcal B_1$. We can show the boundary condition $\frac{\partial B_1(x)}{\partial x}(f(x)+g(x)u(x))+\lambda_1B_1(x)-\epsilon_2\ge 0$ holds by exploiting the SOS decomposition $\frac{\partial B_1(x)}{\partial x}(f(x)+g(x)u(x))+\lambda_1B_1(x)-\epsilon_2=Z(x)^\top QZ(x)$, where $Z(x)=[1,x_1,x_2,x_1x_2,x_1^2,x_2^2]^\top$ and $Q\succeq 0$.

Figure \ref{fig:linsys} shows the zero level set of the quadratic CBC $B_1(x)$. With controller $u_1(x)$ and $u_2(x)$, vector field in \eqref{eq:linsys} guarantees safety with avoiding the unsafe set. For this case, the system admits an ellipsoidal control invariant set. The level sets of $u_1(x)$ and $u_2(x)$ are shown in Figure \ref{fig:SolU1}-\ref{fig:SolU2}. It can be seen that $u(x)\in\mathcal{U}$ for any $x\in\mathcal{B}_1$.

\subsection{Comparison with Control Barrier Functions}
We end this section by a brief comparison between CBF and CBC.

From the point of view of set invariance, the zero-super level set of both CBC and CBF are control invariant. CBC, which is a direct interpretation of control invariance to ensure safety, takes initial conditions into consideration as well - without initial conditions, the CBC formulation is equivalent to CBF. Although the definition of CBF involves the existence of a class-$\mathcal{K}$ function, this, however is a straightforward property that holds for both CBC and CBF.

From the aspect of controller design, the CBF-QP approach relies on a given safe control invariant set, which is free for our approach \eqref{eq:CBCQP2}. For the case where the control invariant set is constructed \emph{a priori}, although the CBF approach endows Lipschitz continuity for the resulting controller, it also introduces unnecessary conservativeness since $\dot B_2(x)$ is bounded by a \emph{fixed} additional relaxation term. Although there are existing works propose to tune the relaxation coefficient $\alpha$ online \cite{xiao2021adaptive}, additional computational complexity and necessary cost trade-off are also introduced. Our approach \eqref{eq:CBCQP1}, on the other hand, is less restricted with an adapted relaxation coefficient $\lambda_1$. For systems with mode switching such as power systems, formulation \eqref{eq:CBCQP1} ensures safety. For continuous controller synthesis, we can also formulate a QP with using $\lambda_1B_1(x)$ as a relaxation term
\begin{equation}\label{eq:CBCQP}
    \begin{split}
        &\mathop{\min }\limits_{u\in\mathcal{U}}~||u - {u^*(x)}||\\
    \mathrm{s.t.}~&\frac{\partial B_1(x)}{\partial x}(f(x)+g(x)u)+\lambda_1B_1(x)\ge 0,  
    \end{split}
\end{equation}
we recall here $\lambda_1$ is a polynomial of $x$, the argument is dropped for simplicity.

\section{Conclusion}\label{sec:conclusion}
In this paper we investigate the problem of safety verification and controller design for safety critical systems. Our approach depends on the evaluation of a control invariant set which encloses the initial set whereas avoiding the unsafe set. We prove that the existence of a control invariant set inside the safe region is sufficient for safety of nonlinear control systems. The formulation only imposes boundary conditions, thus alleviating conservatism. For polynomial systems with semi-algebraic initial and safe sets, we propose an iterative procedure with using SOS program to synthesize the CBC with encoding general affine control limits. We also show that CBC has less conservativeness compared with CBF from numerical simulations. In the future we aim at extending the formulation to discrete time systems.
\bibliographystyle{ieeetr}
\bibliography{ref.bib}

\begin{thebibliography}{10}

\bibitem{guiochet2017safety}
J.~Guiochet, M.~Machin, and H.~Waeselynck, ``Safety-critical advanced robots: A
  survey,'' {\em Robotics and Autonomous Systems}, vol.~94, pp.~43--52, 2017.

\bibitem{lygeros2004reachability}
J.~Lygeros, ``On reachability and minimum cost optimal control,'' {\em
  Automatica}, vol.~40, no.~6, pp.~917--927, 2004.

\bibitem{margellos2011hamilton}
K.~Margellos and J.~Lygeros, ``Hamilton--jacobi formulation for reach--avoid
  differential games,'' {\em IEEE Transactions on automatic control}, vol.~56,
  no.~8, pp.~1849--1861, 2011.

\bibitem{prajna2005necessity}
S.~Prajna and A.~Rantzer, ``On the necessity of barrier certificates,'' {\em
  IFAC Proceedings Volumes}, vol.~38, no.~1, pp.~526--531, 2005.

\bibitem{blanchini1999set}
F.~Blanchini, ``Set invariance in control,'' {\em Automatica}, vol.~35, no.~11,
  pp.~1747--1767, 1999.

\bibitem{prajna2007framework}
S.~Prajna, A.~Jadbabaie, and G.~J. Pappas, ``A framework for worst-case and
  stochastic safety verification using barrier certificates,'' {\em IEEE
  Transactions on Automatic Control}, vol.~52, no.~8, pp.~1415--1428, 2007.

\bibitem{prajna2004safety}
S.~Prajna and A.~Jadbabaie, ``Safety verification of hybrid systems using
  barrier certificates,'' in {\em International Workshop on Hybrid Systems:
  Computation and Control}, pp.~477--492, Springer, 2004.

\bibitem{ames2016control}
A.~D. Ames, X.~Xu, J.~W. Grizzle, and P.~Tabuada, ``Control barrier function
  based quadratic programs for safety critical systems,'' {\em IEEE
  Transactions on Automatic Control}, vol.~62, no.~8, pp.~3861--3876, 2016.

\bibitem{vidyasagar2002nonlinear}
M.~Vidyasagar, {\em Nonlinear systems analysis}.
\newblock SIAM, 2002.

\bibitem{freeman1996control}
R.~A. Freeman and J.~A. Primbs, ``Control lyapunov functions: New ideas from an
  old source,'' in {\em Proceedings of 35th IEEE Conference on Decision and
  Control}, vol.~4, pp.~3926--3931, IEEE, 1996.

\bibitem{xu2017correctness}
X.~Xu, J.~W. Grizzle, P.~Tabuada, and A.~D. Ames, ``Correctness guarantees for
  the composition of lane keeping and adaptive cruise control,'' {\em IEEE
  Transactions on Automation Science and Engineering}, vol.~15, no.~3,
  pp.~1216--1229, 2017.

\bibitem{hsu2015control}
S.-C. Hsu, X.~Xu, and A.~D. Ames, ``Control barrier function based quadratic
  programs with application to bipedal robotic walking,'' in {\em 2015 American
  Control Conference (ACC)}, pp.~4542--4548, IEEE, 2015.

\bibitem{chen2017obstacle}
Y.~Chen, H.~Peng, and J.~Grizzle, ``Obstacle avoidance for low-speed autonomous
  vehicles with barrier function,'' {\em IEEE Transactions on Control Systems
  Technology}, vol.~26, no.~1, pp.~194--206, 2017.

\bibitem{xiao2021adaptive}
W.~Xiao, C.~Belta, and C.~G. Cassandras, ``Adaptive control barrier
  functions,'' {\em IEEE Transactions on Automatic Control}, 2021.

\bibitem{zeng2021safety}
J.~Zeng, B.~Zhang, Z.~Li, and K.~Sreenath, ``Safety-critical control using
  optimal-decay control barrier function with guaranteed point-wise
  feasibility,'' in {\em 2021 American Control Conference (ACC)},
  pp.~3856--3863, IEEE, 2021.

\bibitem{xiao2019control}
W.~Xiao and C.~Belta, ``Control barrier functions for systems with high
  relative degree,'' in {\em 2019 IEEE 58th conference on decision and control
  (CDC)}, pp.~474--479, IEEE, 2019.

\bibitem{chen2021backup}
Y.~Chen, M.~Jankovic, M.~Santillo, and A.~D. Ames, ``Backup control barrier
  functions: Formulation and comparative study,'' {\em arXiv preprint
  arXiv:2104.11332}, 2021.

\bibitem{tan2021high}
X.~Tan, W.~S. Cortez, and D.~V. Dimarogonas, ``High-order barrier functions:
  Robustness, safety and performance-critical control,'' {\em IEEE Transactions
  on Automatic Control}, 2021.

\bibitem{wang2018permissive}
L.~Wang, D.~Han, and M.~Egerstedt, ``Permissive barrier certificates for safe
  stabilization using sum-of-squares,'' in {\em 2018 Annual American Control
  Conference (ACC)}, pp.~585--590, IEEE, 2018.

\bibitem{srinivasan2020synthesis}
M.~Srinivasan, A.~Dabholkar, S.~Coogan, and P.~A. Vela, ``Synthesis of control
  barrier functions using a supervised machine learning approach,'' in {\em
  2020 IEEE/RSJ International Conference on Intelligent Robots and Systems
  (IROS)}, pp.~7139--7145, IEEE, 2020.

\bibitem{robey2020learning}
A.~Robey, H.~Hu, L.~Lindemann, H.~Zhang, D.~V. Dimarogonas, S.~Tu, and
  N.~Matni, ``Learning control barrier functions from expert demonstrations,''
  in {\em 2020 59th IEEE Conference on Decision and Control (CDC)},
  pp.~3717--3724, IEEE, 2020.

\bibitem{clark2021verification}
A.~Clark, ``Verification and synthesis of control barrier functions,'' {\em
  arXiv preprint arXiv:2104.14001}, 2021.

\bibitem{nagumo1942uber}
M.~Nagumo, ``Uber die lage der integralkurven gewokhnlicher di!
  erentialgleichungen,'' {\em Proceedings of the Physico-Mathematical Society
  of Japan}, vol.~24, no.~272, p.~559, 1942.

\bibitem{cunis2021viability}
T.~Cunis and I.~Kolmanovsky, ``Viability, viscosity, and storage functions in
  model-predictive control with terminal constraints,'' {\em Automatica},
  vol.~131, p.~109748, 2021.

\bibitem{prajna2002introducing}
S.~Prajna, A.~Papachristodoulou, and P.~A. Parrilo, ``Introducing sostools: A
  general purpose sum of squares programming solver,'' in {\em Proceedings of
  the 41st IEEE Conference on Decision and Control, 2002.}, vol.~1,
  pp.~741--746, IEEE, 2002.

\bibitem{sostools}
A.~Papachristodoulou, J.~Anderson, G.~Valmorbida, S.~Prajna, P.~Seiler, and
  P.~A. Parrilo, {\em {SOSTOOLS}: Sum of squares optimization toolbox for
  {MATLAB}}.
\newblock \texttt{http://arxiv.org/abs/1310.4716}, 2013.
\newblock Available from \texttt{http://www.eng.ox.ac.uk/control/sostools},
  \texttt{http://www.cds.caltech.edu/sostools} and
  \texttt{http://www.mit.edu/\~{}parrilo/sostools}.

\bibitem{doi:10.1080/10556789908805766}
J.~F. Sturm, ``Using sedumi 1.02, a {MATLAB} toolbox for optimization over
  symmetric cones,'' {\em Optimization Methods and Software}, vol.~11, no.~1-4,
  pp.~625--653, 1999.

\end{thebibliography}
\end{document}